\documentclass[12pt,a4paper]{article}
\usepackage[utf8]{inputenc}
\usepackage{amsmath}
\usepackage{amsfonts}
\usepackage{amssymb}
\usepackage{amsthm}
\usepackage{authblk}
\usepackage{inputenc}
\usepackage{longtable}

\allowdisplaybreaks

\begin{document}
\title{On zeros of the Riemann zeta function}
\author{Xiaolong Wu}
\affil{Ex. Institute of Mathematics, Chinese Academy of Sciences}
\affil{xwu622@comcast.net}
\date{January 20, 2021}
\maketitle

\begin{abstract}

    This paper shows that, in the critical strip, the Riemann zeta function $\zeta(s)$  have the same set of zeros as $F(s):=\int_0^\infty t^{s-1}(e^t+1)^{-1}dt$, and then discusses the behavior of $F(s)$. 
\end{abstract}
\begin{center}\textbf{ \large Introduction}
\end{center}

    Define $\zeta(s)$ as the Riemann zeta function, and define
\begin{equation}\tag{1}
F(s):=\int_{0}^{\infty} \frac{t^{s-1}}{e^t+1} dt,\quad s\in \mathbb{C}.
\end{equation}
    Define the open critical strip by
\begin{equation}\tag{2}
\mathcal{S}:=\{s\in\mathbb{C}:0<\mathcal{R}(s)<1\}
\end{equation}
where $\mathcal{R}(s)$ is the real part of s.

Theorem 1 shows that $\zeta(s)$ and $F(s)$ have the same set of zeros in $\mathcal{S}$. So we can study the zeros of $F(s)$ instead of $\zeta(s)$, and $F(s)$ looks simpler. 

    Write $s=a+ib$. It is obvious that
\begin{equation}\tag{3}
F(s)=\int_{0}^{\infty}\frac{t^{a-1}\cos(b\log t )}{e^t+1} dt+i\int_{0}^{\infty} \frac{t^{a-1}  \sin(b \log t )}{e^t+1} dt=:F_1(a,b)+iF_2(a,b).
\end{equation}
So, $F_2(a,b)\neq 0$ means $F(s)\neq 0$. This paper will work on $F_2$. We split the integration range of $F_2(a,b)$ in to two sub-ranges $[0,R]$ and $[R,\infty)$ for some real $R=e^{2K\pi}$ with integer $K\geq 0, 1 \leq R<\pi$. 

    Define
\begin{equation}\tag{4}
f(t,a):=\frac{t^{a-1}}{e^t+1},\quad g(t):=\frac{1}{e^t+1}.
\end{equation}

Theorem 2 shows that
\begin{equation}\tag{5}
\int_{0}^{R} f(t,a)\sin(b\log t )dt=-\sum_{n=0}^\infty \frac{g^{(n)}(0)b}{n!((n+a)^2+b^2)}R^{n+a}.
\end{equation}

Theorems 3 and 4 discuss properties of $g^{(n)}(0)$, especially, its relation with the Bernoulli numbers. 

Theorem 5 shows that
\begin{equation}\tag{6}
\int_0^R f(t,a)\sin(b\log t )dt>\frac{R^a}{b(e^R+1)}-\frac{0.47177R^a}{b^2}.
\end{equation}
Theorems 6 to 10 study the behavior of
\begin{equation}\notag
\int_R^\infty f(t,a)\sin(b\log t )dt.
\end{equation}
In particular, Theorem 6 shows
\begin{equation}\notag
\int_R^\infty f(t,a)\sin(b\log t )dt=\sum_{k=K}^\infty \int_{\exp(2k\pi/b)}^{\exp((2k+1)\pi/b}h(t,a,b)\sin(b \log t)dt.
\end{equation}
where
\begin{equation}\tag{7}
h(t,a,b):=\frac{t^{a-1}}{e^t+1}-\frac{t^{a-1}e^{a\pi/b}}{e^{te^{\pi/b}}+1}.
\end{equation}
Theorem 7 shows, for $c=e^{\pi/b}$,
\begin{equation}\tag{8}
\int_R^\infty h(t,a,b)dt=\int_R^{cR}\frac{t^{a-1}}{e^t+1}dt.
\end{equation}
Theorem 10 shows, for some range of t, a, b,
\begin{equation}\tag{9}
h(t,a,b)>0.
\end{equation}

    Write the gamma function as
\begin{equation}\tag{10}
\Gamma(s):=\int_{0}^{\infty} \frac{t^{s-1}}{e^t}  dt.
\end{equation}

\begin{center}\textbf{ \large Main Content}
\end{center}

\noindent {\bfseries Theorem 1.}
\textit{ For  $s\in\mathcal{S}$, $\zeta(s)=0$ if and only if $F(s)=0$.
}

\begin{proof}
It is well known that, ( [Titchmarsh 1986] formula 2.4.1 ),
\begin{equation}\tag{1.1}
\zeta(s)=\frac{G(s)}{\Gamma(s)},
\end{equation}
where
\begin{equation}\tag{1.2}
G(s):=\int_{0}^{\infty}\frac{t^{s-1}}{e^t-1} dt.
\end{equation}

Since $\Gamma(s)$ has no zero nor pole in $\mathcal{S}$, $\zeta(s)$ has the same set of zeros as $G(s)$ in $\mathcal{S}$. It is easy to see, write $u=2t$,
\begin{align*}
G(s)-F(s)&=\int_{0}^{\infty}\frac{t^{s-1}}{e^t-1} dt-\int_{0}^{\infty}\frac{t^{s-1}}{e^t+1} dt\\
&=\int_{0}^{\infty}\frac{2t^{s-1}}{e^{2t}-1} dt
=\int_{0}^{\infty}\frac{2(u/2)^{s-1}}{e^u-1} d\frac{u}{2}\\
&=2^{1-s} \int_{0}^{\infty}\frac{u^{s-1}}{e^u-1} du=2^{1-s} G(s).\tag{1.3}
\end{align*}
This means
\begin{equation}\tag{1.4}
(1-2^{1-s} )G(s)=F(s).
\end{equation}
Since $1-2^{1-s}$ has no zero nor pole in $\mathcal{S}$, $F(s)$ has the same set of zeros as $G(s)$, hence as $\zeta(s)$ in $\mathcal{S}$.
\end{proof}

\noindent {\bfseries Theorem 2.}
\textit{ Let $b\geq 100$ be a real, $K>0$ be an integer, and $R=\exp(2K\pi/b)$ such that $1\leq R<\pi$. Let $f(t,a)$ and $g(t)$ be defined as in (4). Then
\begin{equation}\tag{2.1}
\int_0^R f(t,a) \sin(b \log t )dt=-\sum_{n=0}^\infty \frac{g^{(n)}(0)b}{n!((n+a)^2+b^2)}R^{n+a}.
\end{equation}
}
\begin{proof}
Since $g(t)=(e^t+1)^{-1}$ has a pole at $t=\pi i$, the radius of convergence of series
\begin{equation}\tag{2.2}
g(t)=\sum_{n=0}^\infty \frac{g^{(n)}(0)}{n!}t^n
\end{equation}
is $\pi$. We have
\begin{align*}
\int_0^R \frac{t^{s-1}}{e^t+1}dt&=\int_0^Rt^{s-1}\sum_{n=0}^\infty \frac{g^{(n)}(0)}{n!}t^ndt\\
&=\sum_{n=0}^\infty \frac{g^{(n)}(0)}{n!}\int_0^R t^{n+s-1}dt\\
&=\sum_{n=0}^\infty\frac{g^{(n)}(0)}{n!(n+s)}R^{n+s}\\
&=\sum_{n=0}^\infty\frac{g^{(n)}(0)(n+a-ib)}{n!((n+a)^2+b^2)}R^{n+a}e^{ib\log R}\\
&=\sum_{n=0}^\infty \frac{g^{(n)}(0)(n+a)}{n!((n+a)^2+b^2)}R^{n+a}-i\sum_{n=0}^\infty \frac{g^{(n)}(0)b}{n!((n+a)^2+b^2)}R^{n+a}.\tag{2.3}
\end{align*}
On the other hand, we have
\begin{align*}
\int_0^R \frac{t^{s-1}}{e^t+1}dt&=\int_0^R \frac{t^{a-1+bi}}{e^t+1}dt\\
&=\int_0^R \frac{t^{a-1}e^{ib\log t}}{e^t+1}dt\\
&=\int_0^R \frac{t^{a-1}}{e^t+1}\cos(b\log t)dt+i\int_0^R \frac{t^{a-1}}{e^t+1}\sin(b\log t)dt\tag{2.4}
\end{align*}
Equate the imaginary part of (2.3) and (2.4), we get (2.1).
\end{proof}

\noindent {\bfseries Theorem 3.}
\textit{Let $g(t)=(e^t+1)^{-1}$ as in (4), then
\begin{equation}\tag{3.1}
g^{(n)}(0)=\frac{1-2^{n+1}}{n+1}B_{n+1},
\end{equation}
where $B_n$ is the n-th Bernoulli number.
}
\begin{proof}
The Bernoulli numbers have a generating function:
\begin{equation}\tag{3.2}
\frac{t}{e^t-1}=\sum_{n=0}^\infty \frac{B_n}{n!}t^n.
\end{equation}
On the other hand, we have
\begin{equation}\tag{3.3}
\frac{t}{e^t+1}=\sum_{n=0}^\infty \frac{g^{(n)}(0)}{n!}t^{n+1}.
\end{equation}
Since
\begin{equation}\tag{3.4}
\frac{t}{e^t-1}-\frac{t}{e^t+1}=\frac{2t}{e^{2t}-1}=\sum_{n=0}^\infty \frac{B_n}{n!}(2t)^n,
\end{equation}
we have
\begin{align*}
\sum_{n=0}^\infty \frac{g^{(n)}(0)}{n!}t^{n+1}&=\frac{t}{e^t+1}\\
&=\frac{t}{e^t-1}-\frac{2t}{e^{2t}-1}\\
&=\sum_{n=0}^\infty \frac{B_n}{n!}t^n-\sum_{n=0}^\infty \frac{B_n}{n!}(2t)^n\\
&=\sum_{n=0}^\infty \frac{B_n}{n!}(1-2^n)t^n\tag{3.5}
\end{align*}
Equate coefficients of $t^{n+1}$, we get
\begin{equation}\tag{3.6}
\frac{g^{(n)}(0)}{n!}=\frac{B_{n+1}}{(n+1)!}(1-2^{n+1})
\end{equation}
\begin{equation}\tag{3.7}
g^{(n)}(0)=\frac{1-2^{n+1}}{n+1}B_{n+1}.
\end{equation}
\end{proof}

\noindent {\bfseries Theorem 4.}
\textit{
(1) $g^{(2m)}(0)=0, \quad \forall \,m\geq 1$.\\
(2) $g^{(4m+1)}(0)<0, \quad \forall \,m\geq 0$.\\
(3) $g^{(4m-1)}(0)>0, \quad \forall \,m\geq 1$.\\
(4)
\begin{equation}\tag{4.1}
\frac{g^{(2m-1)}(0)}{(2m-1)!}= (-1)^m (1-2^{-2m})\zeta(2m)\frac{2}{\pi^{2m}},\quad \forall \, m\geq 1.
\end{equation}
(5) 
\begin{equation}\tag{4.2}
-\frac{g^{(4m+1)}(0)}{(4m+1)!}\pi^2<\frac{g^{(4m-1)}(0)}{(4m-1)!}<-1.00013814\frac{g^{(4m+1)}(0)}{(4m+1)!}\pi^2,\quad \forall \, m\geq 2.
\end{equation}
(6)
\begin{equation}\tag{4.3}
\frac{2}{\pi^{2m}}<\frac{\vert g^{(2m-1)}(0)\vert}{(2m-1)!}= (1-2^{-2m})\zeta(2m)\frac{2}{\pi^{2m}},\quad \forall \, m\geq 1.
\end{equation}
(7)
\begin{equation}\tag{4.4}
g(t)=\frac{1}{2}-\frac{1}{2}\tanh \frac{t}{2}.
\end{equation}
}
\begin{proof}
(1), (2) and (3) easily follow from well-known properties of the\\ Bernoulli numbers.\\
(4) It is well known that
\begin{equation}\tag{4.5}
\zeta(2m)=\frac{(-1)^{m+1}(2\pi)^{2m}B_{2m}}{2(2m)!},\quad \forall \, m\geq 1.
\end{equation}
Hence by Theorem 3, we have
\begin{align*}
g^{(2m-1)}(0)&=\frac{1-2^{2m}}{2m}B_{2m}\\
&=\frac{1-2^{2m}}{2m}\cdot \frac{(-1)^{m+1}2(2m)!\zeta(2m)}{(2\pi)^{2m}}\\
&=\frac{(-1)^m2(1-2^{-2m})(2m-1)!\zeta(2m)}{\pi^{2m}},\quad \forall \, m\geq 1.\tag{4.6}
\end{align*}
(5) By equation (4.1) we have
\begin{align*}
\frac{g^{(4m-1)}(0)/(4m-1)!}{-g^{(4m+1)}(0)/(4m+1)!}&=\frac{(1-2^{-4m})\zeta(4m)2/\pi^{4m}}{(1-2^{-4m-2})\zeta(4m+2)2/\pi^{4m+2}}\\
&=\frac{(1-2^{-4m})\zeta(4m)}{(1-2^{-4m-2})\zeta(4m+2)}\pi^2,\quad \forall \, m\geq 1.\tag{4.7}
\end{align*}
The last fraction in (4.7) is a decreasing function in m and approaches 1 when $m\rightarrow\infty$. Hence, we have
\begin{equation}\tag{4.8}
1<\frac{(1-2^{-4m})\zeta(4m)}{(1-2^{-4m-2})\zeta(4m+2)}<\frac{(1-2^{-8})\zeta(8)}{(1-2^{-10})\zeta(10)}<1.00013814, \quad \forall \, m\geq 2.
\end{equation}
(6) By [Ge 2012] Theorem 1.1,
\begin{equation}\tag{4.9}
\frac{2(2m)!}{\pi^{2m}(2^{2m}-1)}<\vert B_{2m}\vert =\zeta(2m) \frac{2(2m)!}{(2\pi)^{2m}}.
\end{equation}
Substitute $B_{2m}=2mg^{(2m-1)} (0)/(1-2^{2m} )$, we get
\begin{equation}\tag{4.10}
\frac{2(2m)!}{\pi^{2m}(2^{2m}-1)}<\left\vert \frac{2mg^{(2m-1)}(0)}{1-2^{2m}}\right\vert =\zeta(2m) \frac{2(2m)!}{(2\pi)^{2m}}.
\end{equation}
\begin{equation}\tag{4.11}
\frac{2}{\pi^{2m}}<\left\vert \frac{g^{(2m-1)}(0)}{(2m-1)!}\right\vert = (1-2^{-2m})\zeta(2m) \frac{2}{\pi^{2m}}.
\end{equation}
(7)
\begin{equation}\tag{4.12}
\frac{1}{e^t+1}-\frac{1}{2}=-\frac{1}{2}\cdot \frac{e^t-1}{e^t+1}=-\frac{1}{2}\tanh\frac{t}{2}.
\end{equation}
\end{proof}

\noindent {\bfseries Theorem 5.}
\textit{Let $K>0$ be the largest integer such that $R=e^{2K\pi/b}\leq 2$. Let $a\leq 0.1$ and $b\geq 100$. Then
\begin{equation}\tag{5.1}
-\sum_{n=0}^\infty \frac{g^{(n)}(0)b}{n!((n+a)^2+b^2)}R^{n+a}>-\frac{R^a}{b(e^R+1)}-\frac{0.47177R^a}{b^3}.
\end{equation}
}
\begin{proof}
By maximality of K, we have
\begin{equation}\tag{5.2}
Re^{2\pi/b}=e^{2K\pi/b}e^{2\pi/b}=e^{2(K+1)/b}>2.
\end{equation}
Hence
\begin{equation}\tag{5.3}
R>2e^{-2\pi/b}>2\left(1-\frac{2\pi}{b}\right)=2-\frac{4\pi}{b}.
\end{equation}
\begin{align*}
&-\sum_{n=0}^\infty \frac{g^{(n)}(0)b}{n!((n+a)^2+b^2)}R^{n+a}+\frac{R^a}{b(e^R+1)}\\
&=-\sum_{n=0}^\infty \frac{g^{(n)}(0)b}{n!((n+a)^2+b^2)}R^{n+a}+\frac{R^a}{b}\sum_{n=0}^\infty \frac{g^{(n)}(0)}{n!}R^n\\
&=\frac{R^a}{b}\sum_{n=0}^\infty \frac{g^{(n)}(0)R^n}{n!}\left(\frac{-b^2}{(n+a)^2+b^2}+1\right)\\
&=\frac{R^a}{b}\sum_{n=0}^\infty \frac{g^{(n)}(0)R^n}{n!}\left(\frac{(n+a)^2}{(n+a)^2+b^2}\right)\\
&=\frac{R^a}{b}\left(\frac{a^2}{2(a^2+b^2}-\frac{(1+a)^2R}{4((1+a)^2+b^2)}+\frac{(3+a)^2R^3}{4((3+a)^2+b^2)}-\frac{(5+a)^2R^5}{4((5+a)^2+b^2)}\right)\\
&\quad +\frac{R^a}{b}\sum_{n=7}^\infty \frac{g^{(n)}(0)R^n}{n!}\left(\frac{(n+a)^2}{(n+a)^2+b^2}\right)\\
&>\frac{R^a}{b}\left(-\frac{0.76667}{b^2}+\sum_{m=2}^\infty c_m\right),\tag{5.4}
\end{align*}
where
\begin{align*}
c_m&:=\frac{g^{(4m-1)}(0)R^{4m-1}}{(4m-1)!}\left(\frac{(4m-1+a)^2}{(4m-1+a)^2+b^2}\right)\\
&+\frac{g^{(4m+1)}(0)R^{4m+1}}{(4m+1)!}\left(\frac{(4m+1+a)^2}{(4m+1+a)^2+b^2}\right).\tag{5.5}
\end{align*}
By Theorem 4, we have $g^{(4m-1)}(0)>0, g^{(4m+1)}(0)$ and
\begin{equation}\tag{5.6}
\frac{g^{(4m-1)}(0)}{(4m-1)!}>-\frac{g^{(4m+1)}(0)\pi^2}{(4m+1)!},
\end{equation}
\begin{equation}\tag{5.7}
-\frac{g^{(4m+1)}(0)}{(4m+1)!}>\frac{2}{\pi^{4m+2}}.
\end{equation}
Hence
\begin{align*}
c_m&>-\frac{g^{(4m+1)}(0)R^{4m-1}}{(4m+1)!}\left(\frac{(4m-1+a)^2}{(4m-1+a)^2+b^2}-\frac{(4m+1+a)^2R^2}{(4m+1+a)^2+b^2}\right)\\
&>\frac{2R^{4m-1}}{\pi^{4m+2}}\left(\frac{(4m+1+a)^2}{(4m+1+a)^2+b^2}\left(\frac{7^2\pi^2}{9^2}-R^2\right)\right)\\
&=\frac{2}{Rb^2}\left(\frac{R}{\pi}\right)^{4m}\left(\frac{7^2}{9^2}-\frac{R^2}{\pi^2}\right)\left(\frac{(4m+1+a)^2}{(4m+1+a)^2b^{-2}+1}\right)\tag{5.8}
\end{align*}
Then
\begin{align*}
\sum_{m=0}^\infty c_m&>\frac{2}{Rb^2}\left(\frac{7^2}{9^2}-\frac{R^2}{\pi^2}\right)\sum_{m=2}^\infty \left(\frac{R}{\pi}\right)^{4m}\left(\frac{(4m+1+a)^2}{(4m+1+a)^2b^{-2}+1}\right)\\
&>\frac{0.1996535370}{b^2}\sum_{m=2}^\infty \left(\frac{R}{\pi}\right)^{4m}\left(\frac{9^2}{9^2b^{-2}+1}\right)\\
&>\frac{16.04199633}{b^2}\sum_{m=2}^\infty \left(\frac{R}{\pi}\right)^{4m}\\
&=\frac{16.04199633}{b^2}\cdot\frac{(R/\pi)^8}{1-(R/pi)^4}\\
&>\frac{0.29490227}{b^2}\tag{5.9}
\end{align*}
Hence
\begin{align*}
&-\sum_{n=0}^\infty \frac{g^{(n)}(0)b}{n!((n+a)^2+b^2)}R^{n+a}+\frac{R^a}{b(e^R+1)}\\
&>\frac{R^a}{b}\left(-\frac{0.76667}{b^2}+\frac{0.29490}{b^2}\right)>-\frac{0.47177R^a}{b^3}.\tag{5.10}
\end{align*}
\end{proof}

\textbf{Remark}. Theorem 5 chooses $R\approx 2$ to make the series converge and make the calculation simple.

\noindent {\bfseries Theorem 6.}
\textit{Let a, b be reals such that $0<a<1,b\geq 100$. Let $K>0$ be the largest integer such that $R=e^{2K\pi/b}\leq 2$. Then
\begin{equation}\tag{6.1}
\int_R^\infty f(t,a)\sin(b \log t)dt=\sum_{k=K}^\infty \int_{\exp(2k\pi/b)}^{\exp((2k+1)\pi/b}h(t,a,b)\sin(b \log t)dt.
\end{equation}
}
\begin{proof}
Write
\begin{equation}\tag{6.2}
\int_R^\infty f(t,a)\sin(b \log t)dt=\sum_{k=K}^\infty I_k,
\end{equation}
where
\begin{equation}\tag{6.3}
I_k:=\int_{\exp(2k\pi/b)}^{\exp((2(k+1)\pi/b)}f(t,a)\sin(b \log t)dt.\end{equation}
Then split the range of $I_k$ in to two sub-intervals according to the sign of $\sin(b \log t)$. That is
\begin{equation}\tag{6.4}
I_k=I_{k,1}+I_{k,2},
\end{equation}
where
\begin{equation}\tag{6.5}
I_{k,1}:=\int_{\exp(2k\pi/b)}^{\exp(((2k+1)\pi/b)}f(t,a)\sin(b \log t)dt,
\end{equation}
\begin{equation}\tag{6.6}
I_{k,2}:=\int_{\exp((2k+1)\pi/b)}^{\exp((2(k+1)\pi/b)}f(t,a)\sin(b \log t)dt.
\end{equation}
By changing variable $u=te^{-\pi/b}$ in $I_{k,2}$, we get
\begin{align*}
I_{k,2}&=\int_{\exp(2k\pi/b)}^{\exp((2k+1)\pi/b)}f(ue^{\pi/b},a)\sin(b \log (ue^{\pi/b}))d(ue^{\pi/b})\\
&=e^{\pi/b}\int_{\exp(2k\pi/b)}^{\exp((2k+1)\pi/b)}f(ue^{\pi/b},a)\sin(b \log u+\pi)du\\
&=-e^{\pi/b}\int_{\exp(2k\pi/b)}^{\exp((2k+1)\pi/b)}f(ue^{\pi/b},a)\sin(b \log u)du.\tag{6.7}
\end{align*}
Hence
\begin{align*}
I_k&=\int_{\exp(2k\pi/b)}^{\exp((2k+1)\pi/b)}(f(t,a)-e^{\pi/b}f(te^{\pi/b},a)\sin(b \log t)dt\\
&=\int_{\exp(2k\pi/b)}^{\exp((2k+1)\pi/b)}h(t,a,b)\sin(b \log t)dt\tag{6.8}
\end{align*}
\end{proof}

\noindent {\bfseries Theorem 7.}
\textit{Let a, b, R be reals such that $0<a<1,b>0,R\geq 1$. Let $c:=e^{\pi/b}$. Then
\begin{equation}\tag{7.1}
\int_R^\infty h(t,a,b)dt=\int_R^{cR}\frac{t^{a-1}}{e^t+1}dt.
\end{equation}
\begin{equation}\tag{7.2}
(c-1)\frac{c^{a-1}R^a}{e^{cR}+1}<\int_R^\infty h(t,a,b)dt<(c-1)\frac{R^a}{e^R+1}.
\end{equation}
}
\begin{proof}
Let $0\leq \alpha<\beta$ be reals, $u=ct, t=c^{-1} u$. Then
\begin{equation}\tag{7.3}
\int_\alpha^\beta\frac{t^{a-1}}{e^{ct}+1}dt=\int_{c\alpha}^{c\beta}\frac{c^{1-a}u^{a-1}}{e^u+1}d(c^{-1}u)=c^{-a}\int_{c\alpha}^{c\beta}\frac{u^{a-1}}{e^u+1}du.
\end{equation}
Hence
\begin{align*}
\int_R^\infty h(t,a,b)dt&=\int_R^\infty\left(\frac{t^{a-1}}{e^t+1}-\frac{t^{a-1}e^{a\pi/b}}{e^{te^{\pi/b}}+1}\right) dt\\
&=\int_R^\infty\frac{t^{a-1}}{e^t+1}dt-c^a\int_R^\infty\frac{t^{a-1}}{e^{ct}+1}dt\\
&=\int_R^\infty\frac{t^{a-1}}{e^t+1}dt-c^a\left(c^{-a}\int_{cR}^\infty\frac{t^{a-1}}{e^t+1}dt\right)\\
&=\int_R^\infty\frac{t^{a-1}}{e^t+1}dt-\int_{cR}^\infty\frac{t^{a-1}}{e^t+1}dt\\
&=\int_R^{cR}\frac{t^{a-1}}{e^t+1}dt\tag{7.4}
\end{align*}
Since $t^{a-1}/(e^t+1)$ is decreasing in t, we have
\begin{equation}\tag{7.5}
\int_R^{cR}\frac{t^{a-1}}{e^t+1}dt>(cR-R)\frac{(cR)^{a-1}}{e^{cR}+1}=(c-1)\frac{c^{a-1}R^a}{e^{cR}+1},
\end{equation}
\begin{equation}\tag{7.6}
\int_R^{cR}\frac{t^{a-1}}{e^t+1}dt>(cR-R)\frac{R^{a-1}}{e^R+1}=(c-1)\frac{R^a}{e^R+1}.
\end{equation}
\end{proof}

\noindent {\bfseries Theorem 8.}
\textit{Let k be an integer, $b\geq 10$ be a real. Define the average
\begin{equation}\tag{8.1}
A:=\frac{1}{e^{(2k+1)\pi/b}-e^{2k\pi/b}}\int_{exp(2k\pi/b)}^{exp((2k+1)\pi/b)}\sin (b \log t)dt.
\end{equation}
Then
\begin{equation}\tag{8.2}
\frac{2}{\pi}-\frac{2}{\pi b^2}<A<\frac{2}{\pi}.
\end{equation}
}
\begin{proof}
Substitute $u=b \log t, t=e^{u/b}$, we get
\begin{align*}
\int_{exp(2k\pi/b)}^{exp((2k+1)\pi/b)}\sin (b \log t)dt&=\frac{1}{b}\int_{2k\pi/b}^{(2k+1)\pi/b}e^{u/b}\sin u du\\
&=\frac{e^{u/b}}{b(b^{-2}+1)}\left[\frac{1}{b}\sin u-\cos u\right]_{2k\pi}^{(2k+1)\pi}\\
&=\frac{e^{(2k+1)\pi/b}+e^{2k\pi/b}}{b(b^{-2}+1)}.\tag{8.3}
\end{align*}
So the average is
\begin{equation}\tag{8.4}
A=\frac{e^{(2k+1)\pi/b}+e^{2k\pi/b}}{b(b^{-2}+1)(e^{(2k+1)\pi/b}-e^{2k\pi/b})}=\frac{1+e^{-\pi/b}}{b(b^{-2}+1)(1-e^{-\pi/b})}.
\end{equation}
Since
\begin{align*}
1+e^{-\pi/b}-\frac{2b}{\pi}(1-e^{-\pi/b})&=1+\sum_{n=0}^\infty\frac{(-\pi)^n}{n!b^n}+\frac{2b}{\pi}\sum_{n=1}^\infty\frac{(-\pi)^n}{n!b^n}\\
&=2+\sum_{n=1}^\infty\frac{(-\pi)^n}{n!b^n}-2+2\sum_{n=1}^\infty\frac{(-1)^{n+1}\pi^n}{(n+1)!b^n}\\
&=\sum_{n=1}^\infty\frac{(-\pi)^n}{b^n}\left(\frac{1}{n!}-\frac{2}{(n+1)!}\right)\\
&=\sum_{n=1}^\infty\frac{(-\pi)^n}{b^n}\cdot\frac{n-1}{(n+1)!}\\
&=\frac{\pi^2}{3!b^2}-\cdots,\tag{8.5}
\end{align*}
we have
\begin{align*}
\frac{2}{\pi}&<\frac{1+e^{-\pi/b}}{b(1-e^{-\pi/b})}\\
&<\frac{2}{\pi}+\frac{\pi^2}{3!b^2b(\pi/b-(\pi/b)^2/2)}\\
&=\frac{2}{\pi}+\frac{\pi}{3!b^2(1-\pi/(2b))}\\
&\leq \frac{2}{\pi}+\frac{\pi}{3!b^2(1-\pi/20)}\\
&<\frac{2}{\pi}+\frac{\pi}{3!b^2(5/6)}=\frac{2}{\pi}+\frac{\pi}{5b^2}, \quad \forall\, b\geq 10.\tag{8.6}
\end{align*}
Hence
\begin{align*}
A&<\frac{1}{b^{-2}+1}\left(\frac{2}{\pi}+\frac{\pi}{5b^2}\right)\\
&<\left(1-\frac{1}{b^2}+\frac{1}{b^4}\right)\left(\frac{2}{\pi}+\frac{\pi}{5b^2}\right)\\
&=\frac{2}{\pi}+\frac{\pi}{5b^2}-\frac{2}{\pi b^2}-\frac{\pi}{5b^4}+\frac{2}{\pi b^4}+\frac{\pi}{5b^6}\\
&=\frac{2}{\pi}-\frac{0.0083}{b^2}+\frac{0.0083}{b^4}+\frac{0.6283}{b^6}<\frac{2}{\pi},\tag{8.7}
\end{align*}
\begin{equation}\tag{8.8}
A>\frac{1}{b^{-2}+1}\left(\frac{2}{\pi}\right)>\frac{2}{\pi}\left(1-\frac{1}{b^2}\right).
\end{equation}
\end{proof}

\noindent {\bfseries Theorem 9.}
\textit{Let k be an integer, $b\geq 100$ be a real. Write $t_0=\exp(2k\pi/b)$, $t_1=\exp((2k+1)\pi/b)$. Let $h(t)$ be a continuous function and
\begin{equation}\tag{9.1}
M_1=\min_{t_0\leq t\leq t_1}h(t),\quad\quad M_2=\max_{t_0\leq t\leq t_1}h(t).
\end{equation}
Then
\begin{equation}\tag{9.2}
\left(\frac{2}{\pi}-\frac{2}{\pi b^2}\right)\int_{t_0}^{t_1}M_1dt<\int_{t_0}^{t_1}h(t)\sin (b \log t)dt<\frac{2}{\pi}\int_{t_0}^{t_1}M_2dt.
\end{equation}
}
\begin{proof}
Since $\sin (b \log t)\geq 0$ in interval $[t_0,t_1]$, we have, by Theorem 8,
\begin{align*}
\int_{t_0}^{t_1}h(t)\sin (b \log t)dt&\geq \int_{t_0}^{t_1}M_1\sin (b \log t)dt\\
&>M_1\left(\frac{2}{\pi}-\frac{2}{\pi b^2}\right)(t_1-t_0)\\
&=\left(\frac{2}{\pi}-\frac{2}{\pi b^2}\right)\int_{t_0}^{t_1}M_1dt.\tag{9.3}
\end{align*}
The other inequality can be proved similarly.
\end{proof}

\noindent {\bfseries Theorem 10.}
\textit{Let t, a, b be reals such that $t\geq 1, 0<a\leq e/(e+1)=0.731\cdots$ and $b>0$. We have
\begin{equation}\tag{10.1}
h(t,a,b)>0.
\end{equation}
}
\begin{proof}
We need only to prove that
\begin{equation}\tag{10.2}
\frac{1}{e^t+1}-\frac{e^{a\pi/b}}{e^{te^{\pi/b}}+1}>0.
\end{equation}
\begin{align*}
\frac{1}{e^t+1}-\frac{e^{a\pi/b}}{e^{te^{\pi/b}}+1}&=\frac{e^{te^{\pi/b}}+1-e^{t+a\pi/b}-e^{a\pi/b}}{(e^t+1)(e^{te^{\pi/b}}+1)}\\
&>\frac{e^{t(1+\pi/b)}+1-e^{t+a\pi/b}-e^{a\pi/b}}{(e^t+1)(e^{te^{\pi/b}}+1)}\\
&=\frac{e^t(e^{t\pi/b}-e^{a\pi/b})-(e^{a\pi/b}-1)}{(e^t+1)(e^{te^{\pi/b}}+1)}\\
&>\frac{e^t(t\pi/b-a\pi/b)-a\pi/b}{(e^t+1)(e^{te^{\pi/b}}+1)}\\
&=\frac{\pi}{b}\cdot\frac{te^t-ae^t-a}{(e^t+1)(e^{te^{\pi/b}}+1)}\tag{10.3}
\end{align*}
Since $te^t-ae^t-a>0$ when $t\geq 1,a\leq e/(e+1)=0.731\cdots$, (10.2) holds.
\end{proof}

\begin{center}
{\bfseries \large References}
\end{center}

\noindent 
{[}Ge 2012{]} Hua-feng Ge. \textit{New Sharp Bounds for the Bernoulli Numbers and Refinement of Becker-Stark Inequalities}. Journal of Applied Mathematics Vol 2012, (2012)\\
{[}Titchmarsh 1986{]} E.C. Titchmarsh and D.R. Heath Brown. \textit{The Theory of the Riemann Zeta-function}. Oxford University, 2nd edition (1986)

\begin{center}
{\bfseries \large Appendix}
\end{center}
$g^{(n)} (0)$ for $n\leq 15$.
\begin{equation}\notag
g^{(0)}(0)=\frac{1}{2},\quad g^{(1)}(0)=-\frac{1}{4}, \quad g^{(3)} (0)=\frac{1}{8},
\end{equation}
\begin{equation}\notag
g^{(5)} (0)=-\frac{1}{4} \quad g^{(7)} (0)=\frac{17}{16}, \quad  g^{(9)} (0)=-\frac{31}{4},
\end{equation}
\begin{equation}\notag
g^{(11) } (0)=\frac{691}{8},  \quad g^{(13)}(0) =-\frac{5461}{4},  \quad g^{(15)}(0) =\frac{929569}{32}.
\end{equation}

\end{document}